\documentclass{amsart}
\usepackage{amsfonts,amssymb,amsmath,amsthm}
\usepackage{url}
\usepackage{enumerate}

\makeatletter
\@namedef{subjclassname@2010}{  \textup{2010} Mathematics Subject Classification.}
\makeatother
\newtheorem{thrm}{Theorem}[section]

\theoremstyle{definition}
\newtheorem{definition}[thrm]{Definition}
\newtheorem{remark}[thrm]{Remark}
\numberwithin{equation}{section}
\email{a.akhlidj@fsac.ac.ma}
\email{samkabbaj@yahoo.fr}
\email{hzoubeir2014@gmail.com}
\input{tcilatex}

\begin{document}
\address{ }
\author{Abdellatif Akhlidj$^{(1)}$}
\author{Samir Kabbaj$^{(2)}$}
\author{Hicham Zoubeir$^{(3)}$}
\address{(1) Hassan II University, Faculty of Sciences, Casablanca, Morocco.%
\\
(2) Ibn Tofail University, Department of Mathematics, Laboratory of
Mathematical Analysis and Noncommutative Geometry,\\
Faculty of Sciences, P. O. B : $133,$ Kenitra, Morocco.\\
(3) Ibn Tofail University, Department of Mathematics,\\
Faculty of Sciences, P. O. B : $133,$ Kenitra, Morocco.}
\title[A generalised G\u{a}vru\c{t}a stability of cohomological equations]{A
generalized G\u{a}vru\c{t}a stability of cohomological equations in
nonquasianalytic Carleman classes}

\begin{abstract}
In this paper we introduce the notion of generalized Gavr\`{u}ta stability
of functional equations in order to study, in the framework of a
nonquasianalytic Carleman class, the stability of a class of cohomological
equations.
\end{abstract}

\dedicatory{$\emph{This}$ $\emph{modest}$ $\emph{work}$ $\emph{is}$ $\emph{%
dedicated}$ $\emph{to}$ $\emph{the}$ $\emph{memory}$ $\emph{of}$ $\emph{our}$
$\emph{beloved}$ $\emph{master}$ $\emph{Ahmed}$ $\emph{Intissar}$ $\emph{%
(1951-2017),}$ $\emph{a}$ $\emph{distinguished}$ $\emph{professor,}$ $\emph{a%
}$ $\emph{brilliant}$ $\emph{mathematician,}$ $\emph{a}$ $\emph{man}$ $\emph{%
with}$ $\emph{a}$ $\emph{golden}$ $\emph{heart.}$\emph{\ }}
\maketitle

\section{ Introduction}

The important concept of stability of a functional equation was first
introduced by Ulam in $1940$ when he asked in a talk before the Mathematics
Club of the University of Wisconsin (\cite{ULA}) the following question :

\textit{"Let }$G_{1}$\textit{\ be a group and let }$\left( G_{2},d\right) $%
\textit{\ be a metric group}. \textit{Given any} $\varepsilon >0$, \textit{%
does there }\text{\textit{exist\ a }}$\delta >0$\text{ \textit{such that if
a function }}$h:G_{1}\rightarrow $\text{\textit{\  \ }}$G_{2}$\text{\textit{\ 
}}\textit{satisfies the inequality }$:$ 
\begin{equation*}
d\left( h(xy\right) ,h(x)h(y))\leq \delta
\end{equation*}%
for all $x,$ $y\in G_{1}$\textit{, then there exists a homomorphism }$%
H:G_{1}\rightarrow $\text{\textit{\  \ }}$G_{2}$ with $:$ 
\begin{equation*}
d\left( H(x),h(x)\right) \leq \varepsilon
\end{equation*}%
\textit{\ for all }$x\in G_{1}?"$

Hyers (\cite{HYE}) was the first to answer partially this question when he
showed in $1941$ the following result :

\textit{"If }$E_{1},$ $E_{2}$\textit{\ are Banach spaces and }$%
f:E_{1}\rightarrow E_{2}$\textit{\ is a mapping which satisfy, for some
constant }$\delta >0$ \textit{and} \textit{for all }$x,y\in E_{1},$\textit{\
the condition }$:$%
\begin{equation*}
\left \Vert f(x+y)-f(x)-f(y)\right \Vert \leq \delta \mathit{\ }
\end{equation*}%
\textit{then there exists a unique mapping }$T:E_{1}\rightarrow E_{2}$%
\textit{\ such that }$:$%
\begin{equation*}
\mathit{\ }T(x+y)=T(x)+T(y)
\end{equation*}%
\textit{\ for all }$x,$ $y\in E_{1}$ \textit{and }$:$\textit{\ }%
\begin{equation*}
\left \Vert f(x)-T(x)\right \Vert \leq \delta
\end{equation*}%
\textit{for all }$x\in E_{1}."$

In $1978$ Rassias (\cite{RAS}) has generalized the result of Hyers in the
following way\ :

"\textit{Let} $f:E_{1}\rightarrow E_{2}$ \textit{be a mapping between Banach
spaces and let }$p<1$\textit{\ be fixed. If f satisfies, the inequality }$:$%
\textit{\ }%
\begin{equation*}
\left \Vert f\left( x+y\right) -f(x)-f(y)\right \Vert \leq \theta \left(
\left \Vert x\right \Vert ^{p}+\left \Vert y\right \Vert ^{p}\right)
\end{equation*}%
\textit{holds for each }$x,$ $y\in E_{1}\left( \text{\textit{resp. all} }x,%
\text{ }y\in E_{1}\backslash \left \{ 0\right \} \right) $\textit{\ and for
some constant }$\theta >0$\textit{\ }. \textit{Then there exists a unique
mapping} $T$ $:E_{1}\rightarrow E_{2}$ \textit{such that }$:$%
\begin{equation*}
\mathit{\ }T(x+y)=T(x)+T(y)
\end{equation*}%
\textit{\ for all }$x,$ $y\in E_{1}$ \textit{and }$:$\textit{\ }%
\begin{equation*}
\left \Vert T\left( x\right) -f\left( x\right) \right \Vert \leq \frac{%
2\theta }{\left \vert 2-2^{p}\right \vert }\left \Vert x\right \Vert ^{p}
\end{equation*}%
\textit{for all} $x\in E_{1}$ (\textit{resp. all }$x\in E_{1}\backslash
\left \{ 0\right \} ).$\textit{\ If in addition, }$t\mapsto f\left(
tx\right) $\textit{\ is continuous for each fixed }$x\in E_{1}$\textit{,
then }$T$\textit{\ is linear.}$"$

In $1994$ Gavruta (\cite{GAV}) has given a new generalization of the
Hyers-Ulam-Rassias stability of approximately additive mappings. In fact he
showed that :

\textit{"Let }$G$\textit{\ be an abelian group and }$(X,\left \Vert .\right
\Vert )$\textit{\ a Banach space. Let }$\varphi :G\times G\rightarrow 
%TCIMACRO{\U{211d} }%
%BeginExpansion
\mathbb{R}
%EndExpansion
^{+}$\textit{\ a mapping satisfying, for all }$x,y\in G,$\textit{\ the
condition }$:$\textit{\ }%
\begin{equation*}
\widetilde{\varphi }(x,y):=\mathit{\  \ }\sum \limits_{k=0}^{+\infty
}2^{-k}\varphi (2^{k}x,2^{k}y)<+\infty
\end{equation*}%
\textit{\ Let }$f:G\rightarrow X$ \textit{be a mapping} \textit{which\
fullfiles, for each }$x,$ $y\in G,$ \textit{the condition }$:$\textit{\ }%
\begin{equation*}
\left \Vert f(x+y)-f(x)-f(y)\right \Vert \leq \varphi (x,y)
\end{equation*}%
\textit{\ Then there exists a unique mapping }$T:G\rightarrow X$\textit{\
such that }$:$\textit{\ }%
\begin{equation*}
T(x+y)=T(x)+T(y)
\end{equation*}%
\textit{\ for all }$x,$ $y\in G$\textit{\ and }$:$\textit{\ }%
\begin{equation*}
\left \Vert f(x)-T(x)\right \Vert \leq \frac{1}{2}\widetilde{\varphi }(x,x)
\end{equation*}%
\textit{\ for all }$x\in G.$\textit{\ }$"$

In this paper we introduce the notion of generalized Gavr\`{u}ta stability
of functional equations in order to study, in the framework of a
nonquasianalytic Carleman class $C_{M}\{%
%TCIMACRO{\U{211d} }%
%BeginExpansion
\mathbb{R}
%EndExpansion
\},$ the stability of the so-called cohomological equation $(E_{\psi ,\chi
}):$%
\begin{equation*}
(E_{\psi ,\chi }):f-(f\circ \psi )=\chi
\end{equation*}%
where $f$ is the unknown function and $\chi :%
%TCIMACRO{\U{211d} }%
%BeginExpansion
\mathbb{R}
%EndExpansion
\rightarrow 
%TCIMACRO{\U{2102} }%
%BeginExpansion
\mathbb{C}
%EndExpansion
$, $\psi :$ $%
%TCIMACRO{\U{211d} }%
%BeginExpansion
\mathbb{R}
%EndExpansion
\rightarrow 
%TCIMACRO{\U{211d} }%
%BeginExpansion
\mathbb{R}
%EndExpansion
$ are a given functions belonging to $C_{M}\{%
%TCIMACRO{\U{211d} }%
%BeginExpansion
\mathbb{R}
%EndExpansion
\}.$ Let us recall that cohomological equations play a fundamental role in
the study of dynamical systems. Indeed, the study of certain forms of
invariance, rigidity and stability of dynamical systems can be reduced to
the investigation of the solvability in certain regularity classes of some
cohomological equations (\cite{ALE}, \cite{COR}-\cite{DEL}, \cite{FOR}-\cite%
{LIV}, \cite{MAR1}-\cite{WIL}). However, despite the great interest devoted
to these functional equations, there is at our knowledge a lack of works on
their solvability and their stability in the setting of Carleman classes.
Finally let us pointwise that we were mainly motivated in the preparation of
this paper, by the works (\cite{BEL1}, \cite{BEL2}) of G. Belitskii, E. M.
Dyn'kin and V. Tkachenko. Finally to illustrate our main result, we will
consider the cohomological equations of the form :%
\begin{equation*}
(E_{\chi }):f(x)-f(x+1)=\chi \left( x\right)
\end{equation*}%
which are a particular case of a functional equations called traditionally
difference equations. Let us also recall that, such a functional equations
were studied by numerous authors (\cite{EUL}, \cite{POI}, \cite{GUI}, \cite%
{GUL}, \cite{PIC}, \cite{NOR}, \cite{LAN}, \cite{EVG}, \cite{CHU}, \cite{HAR}%
, \cite{LI}, \cite{TRE}, \cite{PIT}, \cite{JAG}, \cite{AYC}, .etc.) because
of their great importance in applied and fundamental sciences.

\section{Preliminary notes and statement of the main result}

\subsection{Basic notations and main definitions}

For all $x\in 
%TCIMACRO{\U{211d} }%
%BeginExpansion
\mathbb{R}
%EndExpansion
$ we set :%
\begin{equation*}
\left \{ 
\begin{array}{c}
\left \lfloor x\right \rfloor :=\max \left( \{p\in 
%TCIMACRO{\U{2124} }%
%BeginExpansion
\mathbb{Z}
%EndExpansion
:p\leq x\} \right) \\ 
\{x\}:=x-\left \lfloor x\right \rfloor \\ 
\left \lceil x\right \rceil :=\min \left( \{p\in 
%TCIMACRO{\U{2124} }%
%BeginExpansion
\mathbb{Z}
%EndExpansion
:x\leq p\} \right) \\ 
x^{+}:=\max \left( x,0\right)%
\end{array}%
\right.
\end{equation*}

We denote by $B_{r}$, for each $r\in 
%TCIMACRO{\U{2115} }%
%BeginExpansion
\mathbb{N}
%EndExpansion
,$ the Bernouilli number of order $r$ (\cite{DIE}, page $297,299$).

Let $f:S\longrightarrow 
%TCIMACRO{\U{2102} }%
%BeginExpansion
\mathbb{C}
%EndExpansion
$ \ be a function. $\parallel f\parallel _{\infty ,S}$ denotes the quantity :%
\begin{equation*}
||f||_{_{\infty ,S}}:=\sup_{z\in S}|f(s)|
\end{equation*}

Let $X$ be a nonempty set and $F:X\rightarrow X$ a mapping. We denote by $%
F^{\left \langle n\right \rangle }$ for each $n\in 
%TCIMACRO{\U{2115} }%
%BeginExpansion
\mathbb{N}
%EndExpansion
$ the iterate of order $n$ of the mapping $F.$ If $F$ is a bijection, then
we will denote by $F^{\left \langle -1\right \rangle }$the compositional
inverse of the mapping $F$ and for each $n\in 
%TCIMACRO{\U{2115} }%
%BeginExpansion
\mathbb{N}
%EndExpansion
,$ by $F^{\left \langle -n\right \rangle }$ the iterate of order $n$ of the
mapping $F^{\left \langle -1\right \rangle }$.

\begin{definition}
\textit{Let }$E$\textit{\ be a nonempty set, }$F$\textit{\ a nonempty subset
of the set of mappings from }$E$\textit{\ to a metric space }$(V,d),$\textit{%
\ }$\Phi :F\rightarrow F$\textit{\ a given mapping and }$g$\textit{\ a given
element of }$F.$\textit{\ We say that the functional equation }$:$\textit{\ }%
\begin{equation}
\Phi \mathcal{(}y)=g  \label{functequatndefstability}
\end{equation}%
\textit{has the generalized G\u{a}vru\c{t}a stability (GGS)\ in }$F$\textit{%
\ if the following condition is fullfiled }$:$

\textit{For every mapping }$\delta :E\rightarrow 
%TCIMACRO{\U{211d} }%
%BeginExpansion
\mathbb{R}
%EndExpansion
^{+}$\textit{\ there exists a mapping }$\mu :E\rightarrow 
%TCIMACRO{\U{211d} }%
%BeginExpansion
\mathbb{R}
%EndExpansion
^{+}$\textit{\ depending only on }$\delta $\textit{\ and }$\Phi $\textit{\
such that for each mapping }$y\in F$\textit{\ satisfying the inequality }$:$%
\begin{equation*}
d(\Phi \mathcal{(}y)(x),g(x))\leq \delta (x),\text{ }x\in E
\end{equation*}%
\textit{there exists a solution }$z\in F$\textit{\ of the functional
equation (\ref{functequatndefstability}) such that the following condition
holds }$:$%
\begin{equation*}
d(y(x),z(x))\leq \mu (x),\text{ }x\in E
\end{equation*}
\end{definition}

\begin{definition}
\textit{Let }$A:=(A_{n})_{n\text{ }\geq 0\text{ }}$\textit{be a sequence of
strictly positive real numbers.}

\textit{i. The Carleman class }$C_{A}\{%
%TCIMACRO{\U{211d} }%
%BeginExpansion
\mathbb{R}
%EndExpansion
$\textit{\ }$\}$\textit{\ is then the set of all functions }$f:%
%TCIMACRO{\U{211d} }%
%BeginExpansion
\mathbb{R}
%EndExpansion
\rightarrow 
%TCIMACRO{\U{2102} }%
%BeginExpansion
\mathbb{C}
%EndExpansion
$\textit{\ of class }$C^{\infty }$\textit{\ such that }$:$%
\begin{equation*}
||f^{\text{ }(n)}||_{\infty ,K}\leq C_{K}\rho _{K}^{n}A_{n},\text{ }n\in 
%TCIMACRO{\U{2115}}%
%BeginExpansion
\mathbb{N}%
%EndExpansion
\end{equation*}%
\textit{for every compact interval }$K$\textit{\ of }$%
%TCIMACRO{\U{211d} }%
%BeginExpansion
\mathbb{R}
%EndExpansion
$\textit{\ with some constants }$C_{K},\rho _{K}>0.$

ii.\textit{\ The Carleman class }$C_{A}$\textit{\ }$\{%
%TCIMACRO{\U{211d} }%
%BeginExpansion
\mathbb{R}
%EndExpansion
\}$\textit{\ is said to be nonquasinalytic if there exists a nonidentically
vanishing\ function }$f_{0}\in C_{A}$\textit{\ }$\{%
%TCIMACRO{\U{211d} }%
%BeginExpansion
\mathbb{R}
%EndExpansion
\}$\textit{\ such that }$:$%
\begin{equation*}
f_{0}^{\text{ }(n)}(x_{0})=0,\text{ }n\in 
%TCIMACRO{\U{2115}}%
%BeginExpansion
\mathbb{N}%
%EndExpansion
\end{equation*}%
\textit{for some }$x_{0}\in 
%TCIMACRO{\U{211d} }%
%BeginExpansion
\mathbb{R}
%EndExpansion
$\textit{.}

\textit{iii. The sequence }$A$\textit{\ is said to be almost increasing if
there exists a constant }$C>0$\textit{\ such that }$:$%
\begin{equation*}
A_{p}\leq CA_{q}\text{ if }p\leq q
\end{equation*}
\end{definition}

\subsection{Assumptions and related notations}

Along this paper we make the following assumptions :

\begin{itemize}
\item $a<b$ are a fixed real numbers.

\item $M:=(M_{n})_{n\text{ }\geq 0\text{ }}$is a fixed sequence of strictly
positive real numbers{} such that $C_{M}\{%
%TCIMACRO{\U{211d} }%
%BeginExpansion
\mathbb{R}
%EndExpansion
$ $\}$ is nonquasianalytic and the following conditions hold :%
\begin{equation}
\text{the sequence}\left( \left( \frac{M_{n}}{n!}\right) _{\text{ }}^{\frac{1%
}{n}}\right) _{n\in 
%TCIMACRO{\U{2115} }%
%BeginExpansion
\mathbb{N}
%EndExpansion
^{\ast }}\text{ is almost increasing}  \label{seq1}
\end{equation}%
\begin{equation}
1=M_{0}\leq M_{1}  \label{seq2}
\end{equation}%
\begin{equation}
\left( \frac{M_{n+1}}{(n+1)!}\right) ^{2}\leq \frac{M_{n}}{n!}\frac{M_{n+2}}{%
(n+2)!},n\in 
%TCIMACRO{\U{2115} }%
%BeginExpansion
\mathbb{N}
%EndExpansion
\label{seq3}
\end{equation}%
\begin{equation}
\underset{n\in 
%TCIMACRO{\U{2115} }%
%BeginExpansion
\mathbb{N}
%EndExpansion
}{\sup }\left( \frac{M_{n+1}}{(n+1)M_{n}}\right) ^{\frac{1}{n}}<+\infty
\label{seq4}
\end{equation}%
\begin{equation}
\underset{n\rightarrow +\infty }{\lim \inf }\left( \frac{M_{n}}{n!}\right) ^{%
\frac{1}{n}}>0  \label{seq5}
\end{equation}

\item $\psi ,\chi :%
%TCIMACRO{\U{211d} }%
%BeginExpansion
\mathbb{R}
%EndExpansion
\rightarrow 
%TCIMACRO{\U{211d} }%
%BeginExpansion
\mathbb{R}
%EndExpansion
$\ are a fixed functions belonging to the Carleman class $C_{M}\{%
%TCIMACRO{\U{211d} }%
%BeginExpansion
\mathbb{R}
%EndExpansion
\}$ such that the following conditions hold :%
\begin{equation}
\psi (s)>s,\text{ }s\in 
%TCIMACRO{\U{211d} }%
%BeginExpansion
\mathbb{R}
%EndExpansion
\label{psi1}
\end{equation}%
\begin{equation}
\underset{t\rightarrow -\infty }{\lim }\psi (t)=-\infty  \label{psi2}
\end{equation}%
\begin{equation}
\psi ^{\prime }(s)>0,s\in 
%TCIMACRO{\U{211d} }%
%BeginExpansion
\mathbb{R}
%EndExpansion
\label{psi3}
\end{equation}
\end{itemize}

\begin{remark}
It follows from the assumptions (\ref{psi1}-\ref{psi3}) that $\psi $ is a
diffeomorphism from $%
%TCIMACRO{\U{211d} }%
%BeginExpansion
\mathbb{R}
%EndExpansion
$ onto $%
%TCIMACRO{\U{211d} }%
%BeginExpansion
\mathbb{R}
%EndExpansion
$ and that the following relations hold for each $s\in 
%TCIMACRO{\U{211d} }%
%BeginExpansion
\mathbb{R}
%EndExpansion
:$ 
\begin{equation*}
\left \{ 
\begin{array}{c}
\underset{n\rightarrow +\infty }{\lim }\psi ^{\left \langle n\right \rangle
}(s)=+\infty \text{ } \\ 
\underset{n\rightarrow +\infty }{\lim }\psi ^{\left \langle -n\right \rangle
}(s)=-\infty%
\end{array}%
\right.
\end{equation*}
\end{remark}

Let us set for every $s,t\in 
%TCIMACRO{\U{211d} }%
%BeginExpansion
\mathbb{R}
%EndExpansion
:$%
\begin{eqnarray*}
\mathcal{N}_{\psi ,t}^{+}(s) &:&=\min \left( \{p\in 
%TCIMACRO{\U{2115} }%
%BeginExpansion
\mathbb{N}
%EndExpansion
:\psi ^{\left \langle p\right \rangle }(s)\geq t\} \right) \\
\mathcal{N}_{\psi ,t}^{-}(s) &:&=\min \left( \{p\in 
%TCIMACRO{\U{2115} }%
%BeginExpansion
\mathbb{N}
%EndExpansion
:\psi ^{\left \langle -p\right \rangle }(s)\leq t\} \right)
\end{eqnarray*}%
It is clear that the integer valued function : 
\begin{equation*}
\begin{array}{llll}
\mathcal{N}_{\psi ,t}^{+}: & 
%TCIMACRO{\U{211d} }%
%BeginExpansion
\mathbb{R}
%EndExpansion
& \longrightarrow & 
%TCIMACRO{\U{211d} }%
%BeginExpansion
\mathbb{R}
%EndExpansion
^{+} \\ 
& s & \longmapsto & \mathcal{N}_{t}^{+}(s)%
\end{array}%
\end{equation*}%
is increasing while the integer valued function : 
\begin{equation*}
\begin{array}{llll}
\mathcal{N}_{\psi ,t}^{-}: & 
%TCIMACRO{\U{211d} }%
%BeginExpansion
\mathbb{R}
%EndExpansion
& \longrightarrow & 
%TCIMACRO{\U{211d} }%
%BeginExpansion
\mathbb{R}
%EndExpansion
^{+} \\ 
& s & \longmapsto & \mathcal{N}_{t}^{-}(s)%
\end{array}%
\end{equation*}%
is decreasing. We can easily prove that the following inequality holds for
each $s\in 
%TCIMACRO{\U{211d} }%
%BeginExpansion
\mathbb{R}
%EndExpansion
$ and all real numbers $t_{1}<t_{2}$ : 
\begin{equation*}
-\mathcal{N}_{\psi ,t_{1}}^{-}(s)\leq \mathcal{N}_{t_{2}}^{+}(s)-1
\end{equation*}

\subsection{Statement of the main result}

Our main result in this paper is the following theorem.

\begin{theorem}
\end{theorem} \textit{The cohomological equation }$:$\textit{\ }%
\begin{equation*}
(E_{\psi ,\chi })\quad f-(f\circ \psi )=\chi
\end{equation*}%
\textit{has, under the above assumptions, the GGS\ in the Carleman class }$%
C_{M}\{%
%TCIMACRO{\U{211d} }%
%BeginExpansion
\mathbb{R}
%EndExpansion
\}.$\textit{\ More precisely if a function }$y\in C_{M}\{%
%TCIMACRO{\U{211d} }%
%BeginExpansion
\mathbb{R}
%EndExpansion
\}$ \textit{satisfies the condition }$:$\textit{\ }%
\begin{equation*}
|y(s)-y(\psi (s))-\chi (s)|\leq \delta (s),\text{ }s\in 
%TCIMACRO{\U{211d}}%
%BeginExpansion
\mathbb{R}%
%EndExpansion
\end{equation*}%
\textit{where }$\delta :%
%TCIMACRO{\U{211d} }%
%BeginExpansion
\mathbb{R}
%EndExpansion
\rightarrow 
%TCIMACRO{\U{211d} }%
%BeginExpansion
\mathbb{R}
%EndExpansion
^{+},$ \textit{then there exists a solution }$z\in C_{M}\{%
%TCIMACRO{\U{211d} }%
%BeginExpansion
\mathbb{R}
%EndExpansion
\}$\textit{\ of the CE }$(E_{\psi ,\chi })$\textit{\ such that }$:$%
\begin{equation*}
|y(s)-z(s)|\leq \underset{-\mathcal{N}_{\psi ,a}^{-}(s)\leq n\leq \mathcal{N}%
_{\psi ,b}^{+}(s)-1}{\sum }\delta (\psi ^{\left \langle n\right \rangle
}(s)),\text{ }s\in 
%TCIMACRO{\U{211d}}%
%BeginExpansion
\mathbb{R}%
%EndExpansion
\end{equation*}

\section{Proof of the main result}

\subsection{A key result}

We prove first the following result.

\begin{proposition}
\end{proposition} \textit{The cohomological equation }$:$\textit{\ }%
\begin{equation*}
(E_{\psi ,\chi }):\quad f-(f\circ \psi )=\chi
\end{equation*}%
\textit{has a solution }$g$\textit{\ in the Carleman class }$C_{M}\{%
%TCIMACRO{\U{211d} }%
%BeginExpansion
\mathbb{R}
%EndExpansion
\}$ \textit{such that the following inequality holds for every }$s\in 
%TCIMACRO{\U{211d} }%
%BeginExpansion
\mathbb{R}
%EndExpansion
:$%
\begin{equation*}
|g(s)|\leq \underset{-\mathcal{N}_{\psi ,a}^{-}(s)\leq n\leq \mathcal{N}%
_{\psi ,b}^{+}(s)-1}{\sum }\left \vert \chi (\psi ^{\left \langle n\right
\rangle }(s))\right \vert
\end{equation*}

\begin{proof}
Since the Carleman class $C_{M}\{%
%TCIMACRO{\U{211d} }%
%BeginExpansion
\mathbb{R}
%EndExpansion
\}$ is nonquasianalytic there exists, thanks to a result due to S.
Mandelbrojt (\cite{MAN}), a function $\kappa \in C_{M}\{%
%TCIMACRO{\U{211d} }%
%BeginExpansion
\mathbb{R}
%EndExpansion
\}$ such that :%
\begin{equation*}
\left \{ 
\begin{array}{c}
0\leq \kappa (s)\leq 1,\text{ }s\in 
%TCIMACRO{\U{211d} }%
%BeginExpansion
\mathbb{R}
%EndExpansion
\\ 
\kappa (s)=0,\text{ }s\leq a \\ 
\kappa (s)=1,\text{ }s\geq b%
\end{array}%
\right. 
\end{equation*}%
Then let us set for every $s\in 
%TCIMACRO{\U{211d} }%
%BeginExpansion
\mathbb{R}
%EndExpansion
:$%
\begin{equation*}
\chi _{-}(s)=\kappa (s)\chi (s),\text{ }\chi _{+}(s)=(1-\kappa (s))\chi (s)
\end{equation*}%
The functions $\chi _{+}$ and $\chi _{-}$ belong to $C_{M}(%
%TCIMACRO{\U{211d} }%
%BeginExpansion
\mathbb{R}
%EndExpansion
)$ and satisfy the following conditions :%
\begin{equation*}
\left \{ 
\begin{array}{c}
\chi _{+}(s)=0,\text{ }s\geq b \\ 
\chi _{-}(s)=0,\text{ }s\leq a%
\end{array}%
\right. 
\end{equation*}%
Since $\psi $ is a diffeomorphism from $%
%TCIMACRO{\U{211d} }%
%BeginExpansion
\mathbb{R}
%EndExpansion
$ onto $%
%TCIMACRO{\U{211d} }%
%BeginExpansion
\mathbb{R}
%EndExpansion
$ and belongs to the Carleman class $C_{M}\{%
%TCIMACRO{\U{211d} }%
%BeginExpansion
\mathbb{R}
%EndExpansion
\}$ it follows from the assumptions $($\ref{seq1}$)$-$($\ref{seq5}$),$
according to (\cite{ARS}), that $f\circ \psi ,$ $f\circ \psi ^{\left \langle
-1\right \rangle }$ belong to the Carleman class $C_{M}\{%
%TCIMACRO{\U{211d} }%
%BeginExpansion
\mathbb{R}
%EndExpansion
\}$ for each $f\in C_{M}\{%
%TCIMACRO{\U{211d} }%
%BeginExpansion
\mathbb{R}
%EndExpansion
\}.$ Let us then define the operators : 
\begin{equation*}
\left \{ 
\begin{array}{c}
\begin{array}{llll}
\mathcal{L}_{+}: & C_{M}\{%
%TCIMACRO{\U{211d} }%
%BeginExpansion
\mathbb{R}
%EndExpansion
\} & \longrightarrow  & C_{M}\{%
%TCIMACRO{\U{211d} }%
%BeginExpansion
\mathbb{R}
%EndExpansion
\} \\ 
& f & \longmapsto  & f\circ \psi 
\end{array}
\\ 
\begin{array}{llll}
\mathcal{L}_{-}: & C_{M}\{%
%TCIMACRO{\U{211d} }%
%BeginExpansion
\mathbb{R}
%EndExpansion
\} & \longrightarrow  & C_{M}\{%
%TCIMACRO{\U{211d} }%
%BeginExpansion
\mathbb{R}
%EndExpansion
\} \\ 
& f & \longmapsto  & f\circ \psi ^{\left \langle -1\right \rangle }%
\end{array}%
\end{array}%
\right. \text{ }
\end{equation*}%
On the other hand it follows from the assumptions on the function $\psi $
that the sequences\ of intervals $([\psi ^{\left \langle -n\right \rangle
}(b),+\infty \lbrack )_{n\in 
%TCIMACRO{\U{2115} }%
%BeginExpansion
\mathbb{N}
%EndExpansion
}$ and $(]-\infty ,\psi ^{\left \langle n\right \rangle }(a)])_{n\in 
%TCIMACRO{\U{2115} }%
%BeginExpansion
\mathbb{N}
%EndExpansion
}$ are both increasing coverings of $%
%TCIMACRO{\U{211d} }%
%BeginExpansion
\mathbb{R}
%EndExpansion
.$ Thence the following inclusion holds for each compact interval $K:=\left[
\alpha ,\beta \right] $ of $%
%TCIMACRO{\U{211d} }%
%BeginExpansion
\mathbb{R}
%EndExpansion
:$ 
\begin{equation*}
K\subset \lbrack \psi ^{\left \langle -\mathcal{N}_{\psi ,a}^{-}(\alpha
)\right \rangle }(a),+\infty \lbrack \text{ }\cap \text{ }]-\infty ,\psi
^{\left \langle \mathcal{N}_{\psi ,b}^{+}(\beta )\right \rangle }(b)]
\end{equation*}%
Furthermore we have for every $s\in K$ and $n\in 
%TCIMACRO{\U{2115} }%
%BeginExpansion
\mathbb{N}
%EndExpansion
:$%
\begin{equation*}
\left \{ 
\begin{array}{c}
\chi _{+}(\psi ^{\left \langle n\right \rangle }(s))=0\text{ if }n\geq 
\mathcal{N}_{b}^{+}(\beta ) \\ 
((-\chi _{-})\circ \psi ^{\left \langle -1\right \rangle })((\psi
^{\left \langle -1\right \rangle })^{\left \langle n\right \rangle }(s))=0\text{
if }n\geq \mathcal{N}_{a}^{-}(\alpha )%
\end{array}%
\right. 
\end{equation*}%
Thence the series $\sum \mathcal{L}_{+}^{\left \langle n\right \rangle }(\chi
_{+})(s)$ and $\sum \mathcal{L}_{-}^{\left \langle n\right \rangle }(-\chi
_{-}\circ \psi ^{\left \langle -1\right \rangle })(s)$ contain finitely many
non-vanishing terms. Consequently the functions $g_{+},g_{-}:%
%TCIMACRO{\U{211d} }%
%BeginExpansion
\mathbb{R}
%EndExpansion
\rightarrow 
%TCIMACRO{\U{2102} }%
%BeginExpansion
\mathbb{C}
%EndExpansion
$ defined by the relations :%
\begin{equation*}
g_{+}\left( s\right) :=\underset{0\leq n\leq \mathcal{N}_{\psi ,b}^{+}(s)-1}{%
\sum }\mathcal{L}_{+}^{\left \langle n\right \rangle }(\chi _{+})(s)
\end{equation*}%
\begin{equation*}
g_{-}\left( s\right) :=\underset{0\leq n\leq \mathcal{N}_{\psi ,a}^{-}(s)-1}{%
\sum }\mathcal{L}_{-}^{\left \langle n\right \rangle }((-\chi _{-})\circ \psi
^{\left \langle -1\right \rangle })(s)
\end{equation*}%
belong to $C_{M}$\bigskip $\{%
%TCIMACRO{\U{211d} }%
%BeginExpansion
\mathbb{R}
%EndExpansion
\}$ according to (\cite{ARS}). Furthermore easy computations show that the
following estimates hold for each $s\in 
%TCIMACRO{\U{211d} }%
%BeginExpansion
\mathbb{R}
%EndExpansion
$ : 
\begin{equation}
\left \{ 
\begin{array}{c}
\left \vert g_{+}(s)\right \vert \leq \underset{0\leq n\leq \mathcal{N}_{\psi
,b}^{+}(s)-1}{\sum }\left \vert \chi (\psi ^{\left \langle n\right \rangle
}(s))\right \vert  \\ 
\left \vert g_{-}(s)\right \vert \leq \underset{-\mathcal{N}_{\psi
,a}^{-}(s)\leq n\leq -1}{\sum }\left \vert \chi (\psi ^{\left \langle
n\right \rangle }(s))\right \vert 
\end{array}%
\right.   \label{one}
\end{equation}%
It is also clear that we have for every $s\in 
%TCIMACRO{\U{211d} }%
%BeginExpansion
\mathbb{R}
%EndExpansion
:$%
\begin{equation}
\left \{ 
\begin{array}{c}
g_{+}(s)-g_{+}(\psi (s))=\chi _{+}(s) \\ 
g_{-}(s)-g_{-}(\psi (s))=\chi _{-}(s)%
\end{array}%
\right.   \label{two}
\end{equation}%
It follows from (\ref{one}) and (\ref{two}) that the function $g:=g_{+}+g_{-}
$ belongs to $C_{M}$\bigskip $\{%
%TCIMACRO{\U{211d} }%
%BeginExpansion
\mathbb{R}
%EndExpansion
\}$ and is a solution of the cohomological equation $(E_{\psi ,\chi })$ such
that : 
\begin{equation*}
|g(s)|\leq \underset{-\mathcal{N}_{\psi ,a}^{-}(s)\leq n\leq \mathcal{N}%
_{\psi ,b}^{+}(s)-1}{\sum }\left \vert \chi (\psi ^{\left \langle
n\right \rangle }(s))\right \vert ,\text{ }s\in 
%TCIMACRO{\U{211d} }%
%BeginExpansion
\mathbb{R}
%EndExpansion
\end{equation*}

The proof of the proposition is then complete.
\end{proof}

\subsection{End of the proof of the main result}

Let $y\in C_{M}\{%
%TCIMACRO{\U{211d} }%
%BeginExpansion
\mathbb{R}
%EndExpansion
\}$ and $\delta :%
%TCIMACRO{\U{211d} }%
%BeginExpansion
\mathbb{R}
%EndExpansion
\rightarrow 
%TCIMACRO{\U{211d} }%
%BeginExpansion
\mathbb{R}
%EndExpansion
^{+}$. We assume that the following inequality holds for every $s\in 
%TCIMACRO{\U{211d} }%
%BeginExpansion
\mathbb{R}
%EndExpansion
:$%
\begin{equation*}
\left \vert y(s)-y(\psi (s))-\chi (s)\right \vert \leq \delta (s)
\end{equation*}%
Let us then consider the function : 
\begin{equation*}
\begin{array}{llll}
\varphi : & 
%TCIMACRO{\U{211d} }%
%BeginExpansion
\mathbb{R}
%EndExpansion
& \longrightarrow & 
%TCIMACRO{\U{2102} }%
%BeginExpansion
\mathbb{C}
%EndExpansion
\\ 
& s & \longmapsto & y(s)-y(\psi (s))-\chi (s)%
\end{array}%
\end{equation*}%
Then thanks to (\cite{ARS}), the function $\varphi $ belongs to the Carleman
class $C_{M}\{%
%TCIMACRO{\U{211d} }%
%BeginExpansion
\mathbb{R}
%EndExpansion
\}.$ According to the proposition 3, there exists a function $\ h\in $ $%
C_{M}\{%
%TCIMACRO{\U{211d} }%
%BeginExpansion
\mathbb{R}
%EndExpansion
\}$ such that :%
\begin{eqnarray*}
h(s)-h(\psi (s)) &=&\varphi (s),\text{ }s\in 
%TCIMACRO{\U{211d} }%
%BeginExpansion
\mathbb{R}
%EndExpansion
\\
|h(s)| &\leq &\underset{-\mathcal{N}_{\psi ,a}^{-}(s)\leq n\leq \mathcal{N}%
_{\psi ,b}^{+}(s)-1}{\sum }\left \vert \varphi (\psi ^{\left \langle n\right
\rangle }(s))\right \vert ,\text{ }s\in 
%TCIMACRO{\U{211d}}%
%BeginExpansion
\mathbb{R}%
%EndExpansion
\end{eqnarray*}%
Then the function $z:=y-h$ $\in C_{M}\{%
%TCIMACRO{\U{211d} }%
%BeginExpansion
\mathbb{R}
%EndExpansion
\}$ is a solution of the cohomological equation $(E_{\psi ,\chi }).$
Furthermore we have for each $s\in 
%TCIMACRO{\U{211d} }%
%BeginExpansion
\mathbb{R}
%EndExpansion
:$%
\begin{eqnarray*}
|y(s)-z(s)| &=&|h(s)| \\
&\leq &\underset{-\mathcal{N}_{\psi ,a}^{-}(s)\leq n\leq \mathcal{N}_{\psi
,b}^{+}(s)-1}{\sum }\left \vert \varphi (\psi ^{\left \langle n\right
\rangle }(s))\right \vert \\
&\leq &\underset{-\mathcal{N}_{\psi ,a}^{-}(s)\leq n\leq \mathcal{N}_{\psi
,b}^{+}(s)-1}{\sum }\delta (\psi ^{\left \langle n\right \rangle }(s))
\end{eqnarray*}%
It follows that the cohomological equation $(E_{\psi ,\chi })$ has th GGS in
the Carleman class $C_{M}\{%
%TCIMACRO{\U{211d} }%
%BeginExpansion
\mathbb{R}
%EndExpansion
\}.$

We have then achieved the proof of our main result.

$\square $

\section{Example}

The function $\psi _{0}:s\longmapsto s+1$ satisfies the conditions (\ref%
{psi1}-\ref{psi3}). Furtermore the following relations hold for every $s,$ $%
t\in 
%TCIMACRO{\U{211d} }%
%BeginExpansion
\mathbb{R}
%EndExpansion
:$%
\begin{equation*}
\mathcal{N}_{\psi _{0},t}^{+}(s)=\left \lceil (t-s)^{+}\right \rceil ,\text{ 
}\mathcal{N}_{\psi _{0},t}^{-}(s)=\left \lceil (s-t)^{+}\right \rceil
\end{equation*}%
Then, according to the above main result, the cohomological equation :%
\begin{equation*}
(E_{\psi _{0},\chi }):f(s)-f(s+1)=\chi (s)
\end{equation*}%
where $\chi \in C_{M}\{%
%TCIMACRO{\U{211d} }%
%BeginExpansion
\mathbb{R}
%EndExpansion
\}$ is a given function,\ has the GGS in the Carleman class $C_{M}\{%
%TCIMACRO{\U{211d} }%
%BeginExpansion
\mathbb{R}
%EndExpansion
\}.$ More precisely if $y\in C_{M}\{%
%TCIMACRO{\U{211d} }%
%BeginExpansion
\mathbb{R}
%EndExpansion
\}$ and $\delta :%
%TCIMACRO{\U{211d} }%
%BeginExpansion
\mathbb{R}
%EndExpansion
\rightarrow 
%TCIMACRO{\U{211d} }%
%BeginExpansion
\mathbb{R}
%EndExpansion
^{+}$satisfy the following condition%
\begin{equation*}
|y(s)-y(s+1)-\chi (s)|\leq \delta (s),\text{ }s\in 
%TCIMACRO{\U{211d}}%
%BeginExpansion
\mathbb{R}%
%EndExpansion
\end{equation*}%
then there exists a solution $z\in C_{M}\{%
%TCIMACRO{\U{211d} }%
%BeginExpansion
\mathbb{R}
%EndExpansion
\}$\ of the CE $(E_{\psi _{0},\chi })$ \ such that%
\begin{equation}
|y(s)-z(s)|\leq \underset{n=-\left \lceil (s-a)^{+}\right \rceil }{\overset{%
\left \lceil (b-s)^{+}\right \rceil -1}{\sum }}\delta (s+n),\text{ }s\in 
%TCIMACRO{\U{211d} }%
%BeginExpansion
\mathbb{R}
%EndExpansion
\label{Approx}
\end{equation}

1. If the function $\delta $\ is periodic with period $1,$\ then the
estimate (\ref{Approx}) becomes :%
\begin{equation}
|y(s)-z(s)|\leq \left( \left \lceil (b-s)^{+}\right \rceil +\left \lceil
(s-a)^{+}\right \rceil \right) \delta (s)  \label{apprx0}
\end{equation}

2. If the function $\delta $\ is of class $C^{1}$\ on $%
%TCIMACRO{\U{211d} }%
%BeginExpansion
\mathbb{R}
%EndExpansion
$\ then we can improve the estimate (\ref{Approx}) by means of a special
case of the Euler Mac-Laurin formula (\cite{DIE}, page 302-303). Indeed we
have for each $s\in 
%TCIMACRO{\U{211d} }%
%BeginExpansion
\mathbb{R}
%EndExpansion
:$%
\begin{eqnarray*}
&&\underset{n=-\left \lceil (s-a)^{+}\right \rceil }{\overset{\left \lceil
(b-s)^{+}\right \rceil -1}{\sum }}\delta (s+n) \\
&=&\underset{-\left \lceil (s-a)^{+}\right \rceil }{\overset{\left \lceil
(b-s)^{+}\right \rceil -1}{\int }}\delta (s+t)dt+ \\
&&+\frac{\delta \left( s-\left \lceil (s-a)^{+}\right \rceil \right) +\delta
\left( s+\left \lceil (b-s)^{+}\right \rceil -1\right) }{2}+ \\
&&+\underset{-\left \lceil (s-a)^{+}\right \rceil }{\overset{\left \lceil
(b-s)^{+}\right \rceil -1}{\int }}\left( \{t\}-\frac{1}{2}\right) \delta
^{\prime }(s+t)dt
\end{eqnarray*}%
\begin{eqnarray*}
&=&\underset{s-\left \lceil (s-a)^{+}\right \rceil }{\overset{s+\left \lceil
(b-s)^{+}\right \rceil -1}{\int }}\delta (u)du+ \\
&&+\frac{\delta \left( s-\left \lceil (s-a)^{+}\right \rceil \right) +\delta
\left( s+\left \lceil (b-s)^{+}\right \rceil -1\right) }{2}+ \\
&&+\underset{s-\left \lceil (s-a)^{+}\right \rceil \sigma }{\overset{s+\left
\lceil (b-s)^{+}\right \rceil -1}{\int }}\left( \{u-s\}-\frac{1}{2}\right)
\delta ^{\prime }(u)du
\end{eqnarray*}%
\begin{eqnarray*}
&\leq &\left( \left \lceil (b-s)^{+}\right \rceil +\left \lceil
(s-a)^{+}\right \rceil \right) \left \Vert \delta \right \Vert _{\infty ,%
\mathcal{I}_{s}}+ \\
&&+\frac{\left( \left \lceil (b-s)^{+}\right \rceil +\left \lceil
(s-a)^{+}\right \rceil -1\right) }{2}\left \Vert \delta ^{\prime }\right
\Vert _{\infty ,\mathcal{I}_{s}}
\end{eqnarray*}%
where $\mathcal{I}_{s}$ denotes the interval $\left[ s-\left \lceil
(s-a)^{+}\right \rceil ,s+\left \lceil (b-s)^{+}\right \rceil -1\right] .$
Thence the estimate (\ref{Approx}) entails that :%
\begin{eqnarray}
&&|y(s)-z(s)|  \label{apprx1} \\
&\leq &\left( \left \lceil (b-s)^{+}\right \rceil +\left \lceil
(s-a)^{+}\right \rceil \right) \left \Vert \delta \right \Vert _{\infty ,%
\mathcal{I}_{s}}+  \notag \\
&&+\frac{\left( \left \lceil (b-s)^{+}\right \rceil +\left \lceil
(s-a)^{+}\right \rceil -1\right) }{2}\left \Vert \delta ^{\prime }\right
\Vert _{\infty ,\mathcal{I}_{s}}  \notag
\end{eqnarray}

3. If the function $\delta $ is of class $C^{2r+1}$\ on $%
%TCIMACRO{\U{211d} }%
%BeginExpansion
\mathbb{R}
%EndExpansion
$ $(r\in 
%TCIMACRO{\U{2115} }%
%BeginExpansion
\mathbb{N}
%EndExpansion
^{\ast })$\ then we can improve the estimate (\ref{Approx}) by means of a
the general Euler Mac-Laurin formula (\cite{DIE}, page 303-304). Indeed we
have for each $s\in 
%TCIMACRO{\U{211d} }%
%BeginExpansion
\mathbb{R}
%EndExpansion
:$%
\begin{eqnarray*}
|y(s)-z(s)| &\leq &\underset{n=-\left \lceil (s-a)^{+}\right \rceil }{%
\overset{\left \lceil (b-s)^{+}\right \rceil -1}{\sum }}\delta (s+n) \\
&\leq &\underset{-\left \lceil (s-a)^{+}\right \rceil }{\overset{\left
\lceil (b-s)^{+}\right \rceil -1}{\int }}\delta (s+t)dt+ \\
&&+\frac{\delta \left( s-\left \lceil (s-a)^{+}\right \rceil \right) +\delta
\left( s+\left \lceil (b-s)^{+}\right \rceil -1\right) }{2}+ \\
&&+\underset{j=1}{\overset{r}{\sum }}\frac{B_{j}}{(2j)!}\left \vert 
\begin{array}{c}
\delta ^{(2j-1)}\left( s+\left \lceil (b-s)^{+}\right \rceil -1\right) - \\ 
-\delta ^{(2j-1)}\left( s-\left \lceil (s-a)^{+}\right \rceil \right)%
\end{array}%
\right \vert + \\
&&+\frac{\left( r+\frac{1}{2}\right) B_{r}}{\left( 2r+1\right) !}\underset{%
-\left \lceil (s-a)^{+}\right \rceil }{\overset{\left \lceil (b-s)^{+}\right
\rceil -1}{\int }}\left \vert \delta ^{(2r+1)}(s+t)\right \vert dt
\end{eqnarray*}%
\begin{eqnarray*}
&\leq &\left( \left \lceil (b-s)^{+}\right \rceil +\left \lceil
(s-a)^{+}\right \rceil \right) \left \Vert \delta \right \Vert _{\infty ,%
\mathcal{I}_{s}}+ \\
&&+\underset{j=1}{\overset{r}{\sum }}\frac{B_{j}}{(2j)!}\left( \left \lceil
(b-s)^{+}\right \rceil +\left \lceil (s-a)^{+}\right \rceil \right) \left
\Vert \delta ^{(2j)}\right \Vert _{\infty ,\mathcal{I}_{s}}+ \\
&&+\frac{\left( r+\frac{1}{2}\right) B_{r}}{\left( 2r+1\right) !}\left(
\left \lceil (b-s)^{+}\right \rceil +\left \lceil (s-a)^{+}\right \rceil
-1\right) \left \Vert \delta ^{(2r+1)}\right \Vert _{\infty ,\mathcal{I}_{s}}
\end{eqnarray*}%
Finally the estimate (\ref{Approx}$)$ becomes :%
\begin{eqnarray}
&&|y(s)-z(s)|  \label{apprx2} \\
&\leq &\left( \left \lceil (b-s)^{+}\right \rceil +\left \lceil
(s-a)^{+}\right \rceil \right) \left \Vert \delta \right \Vert _{\infty ,%
\mathcal{I}_{s}}+  \notag \\
&&+\underset{j=1}{\overset{r}{\sum }}\frac{B_{j}}{(2j)!}\left( \left \lceil
(b-s)^{+}\right \rceil +\left \lceil (s-a)^{+}\right \rceil \right) \left
\Vert \delta ^{(2j)}\right \Vert _{\infty ,\mathcal{I}_{s}}+  \notag \\
&&+\frac{\left( r+\frac{1}{2}\right) B_{r}}{\left( 2r+1\right) !}\left(
\left \lceil (b-s)^{+}\right \rceil +\left \lceil (s-a)^{+}\right \rceil
-1\right) \left \Vert \delta ^{(2r+1)}\right \Vert _{\infty ,\mathcal{I}_{s}}
\notag
\end{eqnarray}

\bigskip

\bigskip

\end{document}